\documentclass[a4paper,reqno]{amsart}

\usepackage{amsmath}
\usepackage{paralist}
\usepackage{graphics}
\usepackage{epsfig}
\usepackage{amsfonts}
\usepackage{amssymb}
\usepackage{hyperref,cite}
\usepackage{amsthm}
\usepackage{enumerate}
\usepackage[mathscr]{eucal}
\usepackage{eqlist}

\newtheorem{theorem}{Theorem}[section]
\newtheorem{lemma}[theorem]{Lemma}

\def\ifl{\iffalse }

\numberwithin{equation}{section}

\numberwithin{equation}{section}

\theoremstyle{remark}

\begin{document}
\title[The extremal problem for weighted combined energy and the generalization of Nitsche inequality]
{The extremal problem for weighted combined energy and the generalization of Nitsche inequality}

\author{Xiaogao Feng$^{\ast}$}
\address{College of Mathematics and Information, China West Normal University,
       Nanchong 637009,P.R.CHINA}\email{fengxiaogao603@163.com}

\author{Ruyue Tang}
\address{College of Mathematics and Information, China West Normal University,
       Nanchong 637009,P.R.CHINA}\email{3208266846@qq.com}

\author{Ting Peng}
\address{College of Mathematics and Information, China West Normal University,
       Nanchong 637009,P.R.CHINA}\email{2939358133@qq.com}

\thanks{Research supported by the National Natural Science Foundation of China (Grant Nos.11701459 and 12271218).}
\thanks{$^{\ast}$ Corresponding author.}

%\thanks{$^{\ast}$ Corresponding author.}



\begin{abstract}
We consider the existence and uniqueness of a minimizer of the extremal problem for weighted combined energy between two concentric annuli and obtain that the extremal mapping is a certain radial mapping. Meanwhile, this in turn implies a Nitsche type phenomenon and we get a $\frac{1}{|w|^{\lambda}}-$Nitsche type inequality ($\lambda\neq1$).  As an application, on the basis of the relationship between weighted combined energy and weighted combined distortion, we also investigate the extremal problem for weighted combined distortion on annuli. This extends the result obtained by Kalaj in \cite{Ka1}.
\end{abstract}

\maketitle

\section{Introduction}
\quad

For two positive and distinct constants $r$ and $R$, let
\[\mathbb{A}_{1}=\{z:1\leq |z|\leq r\}  \quad  \text{and} \quad \mathbb{A}_{2}=\{w:1\leq |w|\leq R\} \]
be rounded annuli in the complex plane $\mathbb{C}$.
we consider the class $\mathfrak{H}(\mathbb{A}_{1},\mathbb{A}_{2})$ of all orientation preserving homeomorphisms
$h$ from $\mathbb{A}_{1}$ onto $\mathbb{A}_{2}$ keeping orders of the boundaries, namely,
\[|h(z)|=1 \quad \text{for} \quad |z|=1  \quad \text{and} \quad |h(z)|=R \quad \text{for} \quad |z|=r. \eqno{(1.1)}\]

For $z=te^{i\theta}$, the normal and tangential derivatives of $h$ are
\[h_{N}=h_{t} \quad \text{and} \quad h_{T}=\frac{1}{t}h_{\theta}. \eqno{(1.2)}\]
The derivatives $h_{z}$ and $h_{\overline{z}}$ are expressed as
\[h_{z}=\frac{e^{-i\theta}} {2}(h_{N}-ih_{T}) \quad \text{and} \quad h_{\overline{z}}=\frac{e^{i\theta}} {2}(h_{N}+ih_{T}). \eqno{(1.3)}\]
We find the Hilbert-Schmidet norm of the differential matrix
\[|Dh|^{2}=2(|h_{z}|^{2}+|h_{\overline{z}}|^{2})=|h_{N}|^{2}+|h_{T}|^{2}, \eqno{(1.4)} \]
and the Jocobian determinant of $h$
\[J(z,h)=|h_{z}|^{2}-|h_{\overline{z}}|^{2}=\Im(\overline{h_{N}}h_{T}).\eqno{(1.5)}\]

In \cite{IO2} Iwaniec and Oninnen have considered the $n-$harmonic energy between two annuli. When $n=2$, in \cite{IO1} they gave a new proof.  For a similar problem but for non-circular annuli we refer to the paper \cite{IKKO} and its generalization in \cite{Ka2}. In 2020, Kalaj \cite{Ka1} has generalized the harmonic energy to the combined energy. He studied the following extremal problem
\[\inf_{h\in\mathfrak{H}(\mathbb{A}_{1},\mathbb{A}_{2})} \iint_{\mathbb{A}_{1}}\left(a^{2}|h_{N}|^{2}+b^{2}|h_{T}|^{2}\right)dz\eqno{(1.6)}\]
and proved that the minimizer is certain radial mapping (see Theorem 3.1 in \cite{Ka1}). Very recently, Yang, Tang and Feng \cite{YTF} furthered (1.6) to the weighted case and discussed the following the extremal problem
\[\inf_{h\in\mathfrak{H}(\mathbb{A}_{1},\mathbb{A}_{2})} \iint_{\mathbb{A}_{1}}\left(a^{2}|h_{N}|^{2}+b^{2}|h_{T}|^{2}\right)\frac{1}{|h(z)|^{4}}dz. \eqno{(1.7)}\]

   In this note, we mainly investigate the mapping in $\mathfrak{H}(\mathbb{A}_{1},\mathbb{A}_{2})$ which minimizes the following extremal problem
for weighted combined energy
\[\mathbb{E}_{\lambda}[h]=\iint_{\mathbb{A}_{1}}\left(a^{2}|h_{N}|^{2}+b^{2}|h_{T}|^{2}\right)\frac{1}{|h(z)|^{2\lambda}}dz, \eqno{(1.8)}\]
where $a>0$, $b>0$ and $\lambda\in \mathbb{R}$.

The harmonic energy has close relationship to Nitsche phenomena. In 1962, Nitsche \cite{Ni1} conjectured that a necessary and sufficient condtion for
existence of a harmonic homeomorphism between two annuli $\mathbb{A}_{1}$  and $\mathbb{A}_{2}$ is the following inequality
\[r\leq R+\sqrt{R^{2}-1} \quad  \text{or}  \quad  R\geq\frac{1}{2}(r+\frac{1}{r}). \eqno{(1.9)}\]
After various lower bounds for $R$ (\cite{Ly}, \cite{We}, \cite{Ni2}), the Nitsche conjecture was finally solved by Iwaniec, Kovalev and
Onninen \cite{IKO} in 2011.

Recall that, given a Riemannian metric $\rho$ on $\mathbb{A}_{2}$, a $C^{2}$ homeomorphism $h$ is said to be harmonic with
respect to $\rho$ (or $\rho-$harmonic) if
\[h_{z\overline{z}}(z)+\left(\log \rho^{2}\right)_{w}\circ hh_{z}h_{\overline{z}}=0. \]
On the basis of Nitsche conjecture, Kalaj \cite{Ka4} proposed the so-called $\rho-$Nitsche conjecture as follows

\textbf{$\rho-$Nitsche conjecture}   \emph{If there exists a $\rho-$harmonic homeomorphism of the annuli  $\mathbb{A}_{1}$  onto $\mathbb{A}_{2}$, then}
\[r\leq \exp\left(\int_{1}^{R}\frac{\rho(s)}{\sqrt{s^{2}\rho^{2}(s)-\alpha_{0}}}ds\right), \eqno{(1.10)}\]
\emph{where} $\alpha_{0}=\inf_{1\leq s\leq R}(\rho^{2}(s))s^{2}$.

It should be pointed out that the $\rho-$Nitsche conjecture corresponds to the classical Nitsche conjecture when $\rho=1$. In \cite{Ka3} a partial result is proved on the $\rho-$Nitsche conjecture. In particular, when $\rho(w)=|w|^{-2}$, then inequality (1.10) becomes
\[r\leq R+\sqrt{R^{2}-1}. \eqno{(1.11)}\]
According to this property, Feng and Tang \cite{FT} gave a positive answer when $\rho(w)=|w|^{-2}$.
Recently, Kalaj \cite{Ka1} has considered the extremal problem (1.6) and obtained a
Nitsche type inequality
\[R\geq\cosh(\frac{b}{a}\log r)=\frac{1+r^{\frac{2b}{a}}}{2r^{\frac{b}{a}}}.\eqno{(1.12)}\]
Notice that when $a=b$, (1.12) turns out to be the classic Nitsche inequality (1.9).

Very recently, in \cite{YTF} we have considered the extremal problem (1.7) and got $\frac{1}{|w|^{2}}-$Nitsche type inequality
\[r\leq(R+\sqrt{R^{2}-1})^{\frac{a}{b}}.\eqno{(1.13)}\]
Moreover, we obtain the following result.

\textbf{Theorem A} (see Theorem 3.1 in \cite{YTF}) Under condition (1.13) and among all mappings in  $\mathfrak{H}(\mathbb{A}_{1},\mathbb{A}_{2})$ , the
infimum of
\[\iint_{\mathbb{A}_{1}}\left(a^{2}|h_{N}|^{2}+b^{2}|h_{T}|^{2}\right)\frac{1}{|h(z)|^{4}}dz \]
is attained at the radial mapping
\[h_{\ast}(z)=\frac{2(1+\sqrt{1+\frac{\alpha}{b^{2}}})|z|^{\frac{b}{a}}}{(1+\sqrt{1+\frac{\alpha}{b^{2}}})^{2}-|z|^{\frac{2b}{a}}\frac{\alpha}{b^{2}}}
e^{i\theta},\eqno{(1.14)}\]
where $\alpha$ satisfies
\[R=\frac{2r^{\frac{b}{a}}(1+\sqrt{1+\frac{\alpha}{b^{2}}})}{(1+\sqrt{1+\frac{\alpha}{b^{2}}})^{2}-r^{\frac{2b}{a}}\frac{\alpha}{b^{2}}}. \eqno{(1.15)}\]
The minimizer is unique up to a rotation of annuli.

  In this paper, we shall generalize the extremal problems (1.6) and (1.7) to (1.8) and get the following main theorem.
\begin{theorem}
(I) When $\lambda\neq 1$ in (1.8),
under condition
\[r\leq\left(R^{|\lambda-1|}+\sqrt{R^{2|\lambda-1|}-1}\right)^{\frac{1}{|\lambda-1|}\frac{a}{b}} \eqno{(1.16)}\]
 and among all mapping $\mathfrak{H}(\mathbb{A}_{1},\mathbb{A}_{2})$, the infimum of
\[\iint_{\mathbb{A}_{1}}\left(a^{2}|h_{N}|^{2}+b^{2}|h_{T}|^{2}\right)\frac{dz}{|h(z)|^{2\lambda}}\]
is attained at the radial mapping
\[h^{\ast}_{\lambda}(z)=\frac{2^{\frac{1}{\lambda-1}}(1+\sqrt{1+\frac{\alpha}{b^{2}}})^{\frac{1}{\lambda-1}}|z|^{\frac{b}{a}}}
{\left[(1+\sqrt{1+\frac{\alpha}{b^{2}}})^{2}-|z|^{\frac{2(\lambda-1)b}{a}}\frac{\alpha}{b^{2}}\right]^{\frac{1}{\lambda-1}}}
e^{i\theta},\eqno{(1.17)}\]
where $\alpha$ satisfies
\[R=\frac{2^{\frac{1}{\lambda-1}}(1+\sqrt{1+\frac{\alpha}{b^{2}}})^{\frac{1}{\lambda-1}}r^{\frac{b}{a}}}
{\left[(1+\sqrt{1+\frac{\alpha}{b^{2}}})^{2}-r^{\frac{2(\lambda-1)b}{a}}\frac{\alpha}{b^{2}}\right]^{\frac{1}{\lambda-1}}}
e^{i\theta}.\eqno{(1.18)}\]
That is to say,
\[\iint_{\mathbb{A}_{1}}\left(a^{2}|h_{N}|^{2}+b^{2}|h_{T}|^{2}\right)\frac{dz}{|h(z)|^{2\lambda}}
\geq\iint_{\mathbb{A}_{1}}\left(a^{2}|h^{\ast}_{\lambda N}|^{2}+b^{2}|h^{\ast}_{\lambda T}|^{2}\right)\frac{dz}{|h^{\ast}_{\lambda}(z)|^{2\lambda}}. \eqno{(1.19)}\]
The minimizer is unique up to a rotation of annuli.

(II)When $\lambda=1$ in (1.8), the weighted combined energy $\mathbb{E}_{\lambda}[h]$ attains its minimum for a radial mapping
\[h_{1}^{\ast}(z)=|z|^{\frac{\ln R}{\ln r}}e^{i\theta}.\eqno{(1.20)}\]
The minimizer is unique up to a rotation of annuli.
\end{theorem}

\textbf{Remark 1.2}  \quad When $\lambda=0$, (I) in Theorem 1.1 happens to be Theorem 3.1 in \cite{Ka1}.

\textbf{Remark 1.3}  \quad When $\lambda=2$, (I) in Theorem 1.1 happens to be Theorem A (see Theorem 3.1 in \cite{YTF}).

\textbf{Remark 1.4}  \quad When $\lambda=0$,  the inequality (1.16) becomes (1.12).

\textbf{Remark 1.5}  \quad When $\lambda=2$,  the inequality (1.16) becomes (1.13).

We end this introduction section with the organization of the paper. In section 2
we shall investigate the weighted combined energy of radial mapping and construct radial minimizers of weighted combined energy. In
Section 3, we will give the proof of Theorem 1.1. The last section,
By means of the relationship between weighted combined energy and weighted combined distortion, we also investigate the extremal problem for weighted combined distortion.

\section{weighted combined energy of radial mapping}

For mappings of annuli, it is natural to examine the radial mappings.
Assume
\[h(z)=H(t)e^{i\theta} \quad  \text{for} \quad  z=te^{i\theta},\] where $H:[1,r]\rightarrow[1,R]$ is continuous strictly increasing and satisfies the boundary conditions (1.1). We obtain that
\[\mathbb{E}_{\lambda}[h]=2\pi\int_{1}^{r}\frac{a^{2}t^{2}\dot{H}^{2}(t)+b^{2}H^{2}(t)}{tH^{2\lambda}(t)}dt
=2\pi\int_{1}^{r}F(t,H(t),\dot{H}(t))dt. \eqno{(2.1)}\]
Then
$$\frac{\partial F}{\partial H}=-\frac{2\lambda a^{2}t\dot{H}^{2}(t)}{H^{2\lambda+1}(t)}-\frac{(2\lambda-2)b^{2}}{tH^{2\lambda-1}(t)},$$
$$\frac{\partial F}{\partial \dot{H}}=\frac{2a^{2}t\dot{H}(t)}{H^{2\lambda}(t)},$$
and
$$\frac{d}{dt}\left(\frac{\partial F}{\partial \dot{H}}\right)=\frac{2a^{2}\dot{H}(t)}{H^{2\lambda}(t)}+
\frac{2a^{2}t\ddot{H}(t)}{H^{2\lambda}(t)}-\frac{4\lambda a^{2}t\dot{H}^{2}(t)}{H^{2\lambda+1}(t)}.$$
According to the Euler-Lagrange equation
\[\frac{\partial F}{\partial H}=\frac{d}{dt}\left(\frac{\partial F}{\partial \dot{H}}\right),  \eqno{(2.2)}\]
we get
\[a^{2}t^{2}H(t)\ddot{H}(t)+a^{2}tH(t)\dot{H}(t)+(\lambda-1)b^{2}H^{2}(t)=\lambda a^{2}t^{2}\dot{H}^{2}(t). \eqno{(2.3)}\]
Set $t=e^{x}$, $y=H(t)=H(e^{x})$, then $y'=\dot{H}(t)e^{x}$. Furthermore,
\[y''=\ddot{H}(t)e^{2x}+\dot{H}(t)e^{x}.\]
Then
\[\dot{H}(t)=\frac{y'}{e^{x}} \quad \text{and} \quad \ddot{H}(t)=\frac{y''-y'}{e^{2x}}. \eqno{(2.4)}\]
We put (2.4) into (2.3) and obtain
\[a^{2}yy''+(\lambda-1)b^{2}y^{2}=\lambda a^{2}y'^{2}.  \eqno{(2.5)}\]
Moreover, set $y'=\zeta$, then
\[y''=\frac{d\zeta}{dy}\frac{dy}{dx}=\zeta'\zeta.\eqno{(2.6)}\]
According to (2.5) and (2.6), we have
\[a^{2}\zeta \zeta'+(\lambda-1)b^{2}y=\frac{\lambda}{y}a^{2}\zeta^{2}.\eqno{(2.7)}\]
Thus
\[a^{2}\zeta \zeta'-b^{2}y=\frac{\lambda}{y}(a^{2}\zeta^{2}-b^{2}y^{2}).\eqno{(2.8)} \]
Let
\[\omega=a^{2}\zeta^{2}-b^{2}y^{2},\eqno{(2.9)}\]
then $\omega'=2a^{2}\zeta \zeta'-2b^{2}y$. Hence, the equality (2.8) reduces to
\[\frac{\omega'}{\omega}=\frac{2\lambda}{y}.\eqno{(2.10)}\]
We get $\omega=\alpha y^{2\lambda}$. Together with (2.9), we get
\[\alpha y^{2\lambda}=a^{2}\zeta^{2}-b^{2}y^{2},\eqno{(2.11)}\]
and
\[y'= \frac{dy}{dx}=\zeta=\frac{\sqrt{\alpha y^{2\lambda}+b^{2}y^{2}}}{a} \eqno{(2.12)}\]
which also can be written as
\[dx=\frac{a}{\sqrt{\alpha y^{2\lambda}+b^{2}y^{2}}}dy. \eqno{(2.13)}\]
Thus, for $1\leq y\leq R$, we have
\[x(y)=\int_{1}^{y}\frac{a}{\sqrt{\alpha s^{2\lambda}+b^{2}s^{2}}}ds.\eqno{(2.14)}\]
Furthermore,
\[t=e^{x}=\exp\left(\int_{1}^{y}\frac{a}{\sqrt{\alpha s^{2\lambda}+b^{2}s^{2}}}ds\right)
.\eqno{(2.15)}\]
In order to discuss (2.15), we divide $\lambda$ into three cases: $\lambda>1$, $\lambda<1$ and $\lambda=1$.

When $\lambda>1$ in (2.15), we take $\alpha=-\inf_{1\leq s\leq R}(\frac{b^{2}}{s^{2(\lambda-1)}})=-\frac{b^{2}}{R^{2\lambda-2}}$ and have the following inequality
\begin{eqnarray*}
r&\leq&\exp\left[\frac{1}{1-\lambda}\frac{a}{b}\int_{1}^{R}\frac{d(\frac{b}{s^{\lambda-1}})}{\sqrt{\left(\frac{b}{s^{\lambda-1}}\right)^{2}
-\frac{b^{2}}{R^{2\lambda-2}}}}\right] \\
&=&\left(R^{\lambda-1}+\sqrt{R^{2(\lambda-1)}-1}\right)^{\frac{1}{\lambda-1}\frac{a}{b}}, \quad \quad(\frac{1}{|w|^{\lambda}}-\text{Nitsche type inequality})
\end{eqnarray*}
\[\eqno{(2.16)}\]
which happens to be (1.16) when $\lambda>1$.
Notice that when $\lambda=2$, (2.16) becomes (1.13). Next we will construct the radial mapping $h^{\ast}_{\lambda}(z)$ (see (1.17)). For $\alpha\geq-\frac{b^{2}}{R^{2\lambda-2}}$, we obtain
\[\exp\left(\frac{1}{1-\lambda}\frac{a}{b}\int_{1}^{R}\frac{d(\frac{b}{s^{\lambda-1}})}{\sqrt{(\frac{b}{s^{\lambda-1}})^{2}+\alpha}}\right)
=\left[\frac{R^{\lambda-1}(1+\sqrt{1+\frac{\alpha}{b^{2}}})}{1+\sqrt{1+\frac{R^{2\lambda-2}\alpha}{b^{2}}}}\right]^{\frac{1}{\lambda-1}\frac{a}{b}}. \eqno{(2.17)}\]
Set
\[\varphi(\alpha)=\left[\frac{R^{\lambda-1}(1+\sqrt{1+\frac{\alpha}{b^{2}}})}{1+\sqrt{1+\frac{R^{2\lambda-2}\alpha}{b^{2}}}}\right]
^{\frac{1}{\lambda-1}\frac{a}{b}},\eqno{(2.18)}\]
then
\begin{eqnarray*}
\varphi'(\alpha)&=&\frac{1}{\lambda-1}\frac{a}{b}\left[\frac{R^{^{\lambda-1}}(1+\sqrt{1+\frac{\alpha}{b^{2}}})}
{1+\sqrt{1+\frac{R^{2\lambda-2}\alpha}{b^{2}}}}\right]^{\frac{1}{\lambda-1}\frac{a}{b}-1}\\
&&\times\frac{R}{2b^{2}}
\frac{\sqrt{1+\frac{R^{2\lambda-2}\alpha}{b^{2}}}-R^{2\lambda-2}\sqrt{1+\frac{\alpha}{b^{2}}}+1-R^{2\lambda-2}}
{\sqrt{1+\frac{\alpha}{b^{2}}}\sqrt{1+\frac{R^{2\lambda-2}\alpha}{b^{2}}}\left(1+\sqrt{1+\frac{R^{2\lambda-2}\alpha}{b^{2}}}\right)^{2}}.
\end{eqnarray*}
\[\eqno{(2.19)}\]
Since $\alpha\geq-\frac{b^{2}}{R^{2\lambda-2}}$ and $\lambda>1$, we get
\begin{eqnarray*}
&&\sqrt{1+\frac{R^{2(\lambda-1)}\alpha}{b^{2}}}-R^{2(\lambda-1)}\sqrt{1+\frac{\alpha}{b^{2}}} \\
&=&\frac{\left[1-R^{2(\lambda-1)}\right]\left[1+R^{2(\lambda-1)}+\frac{R^{2(\lambda-1)}\alpha}{b^{2}}\right]}{\sqrt{1+\frac{R^{2(\lambda-1)}\alpha}{b^{2}}}
+R^{2(\lambda-1)}\sqrt{1+\frac{\alpha}{b^{2}}}}\\
&<&0.
\end{eqnarray*}
\[\eqno{(2.20)}\]
Combining (2.19) and (2.20), we have $\varphi'(\alpha)<0$. Thus, for $\alpha\geq-\frac{b^{2}}{R^{2\lambda-2}}$, $\varphi(\alpha)$ is
strictly decreasing. By (2.16) and
\begin{eqnarray*}
\varphi(-\frac{b^{2}}{R^{2\lambda-2}})&=& \left[\frac{R^{\lambda-1}(1+\sqrt{1+\frac{-\frac{b^{2}}{R^{2\lambda-2}}}{b^{2}}})}{1+\sqrt{1+\frac{R^{2\lambda-2}}{b^{2}}(-\frac{b^{2}}{R^{2\lambda-2}})}}
\right]^{\frac{1}{\lambda-1}\frac{a}{b}}\\
&=&\left(R^{\lambda-1}+\sqrt{R^{2\lambda-2}-1}\right)^{\frac{1}{\lambda-1}\frac{a}{b}}.
\end{eqnarray*}
\[\eqno{(2.21)}\]
We obtain that there exists a unique $\alpha\geq-\frac{b^{2}}{R^{2\lambda-2}}$ such that $\varphi(\alpha)=r$, that is,
\[\left[\frac{R^{\lambda-1}(1+\sqrt{1+\frac{\alpha}{b^{2}}})}{1+\sqrt{1+\frac{R^{2\lambda-2}\alpha}{b^{2}}}}\right]^{\frac{1}{\lambda-1}\frac{a}{b}}=r. \eqno{(2.22)}\]
We conclude that
\[R=\frac{2^{\frac{1}{\lambda-1}}(1+\sqrt{1+\frac{\alpha}{b^{2}}})^{\frac{1}{\lambda-1}}r^{\frac{b}{a}}}{\left[(1+\sqrt{1+\frac{\alpha}{b^{2}}})^{2}
-r^{\frac{2(\lambda-1)b}{a}}\frac{\alpha}{b^{2}}\right]^{\frac{1}{\lambda-1}}}. \eqno{(2.23)}\]
So we get the radial mapping
\[h^{\ast}_{\lambda}(z)=\frac{2^{\frac{1}{\lambda-1}}(1+\sqrt{1+\frac{\alpha}{b^{2}}})^{\frac{1}{\lambda-1}}|z|^{\frac{b}{a}}}
{\left[(1+\sqrt{1+\frac{\alpha}{b^{2}}})^{2}-|z|^{\frac{2(\lambda-1)b}{a}}\frac{\alpha}{b^{2}}\right]^{\frac{1}{\lambda-1}}}
e^{i\theta},\eqno{(2.24)}\]
where $\alpha$ satisfies (2.23). By (2.23) and direct computation we obtain
\[\mathbb{E}[h^{\ast}_{\lambda}]=\frac{2\pi ab}{\lambda-1}\left(\sqrt{1+\frac{\alpha}{b^{2}}}-\frac{1}{R^{\lambda-1}}\sqrt{\frac{1}{R^{2(\lambda-1)}}+\frac{\alpha}{b^{2}}}\right).\eqno{(2.25)}\]

  When $\lambda<1$ in (2.15),  we take $\alpha=-\inf_{1\leq s\leq R}\left(\frac{b^{2}}{s^{2(\lambda-1)}}\right)=-b^{2}$ and have the following
inequality
\begin{eqnarray*}
r&\leq&\exp\left[\frac{1}{1-\lambda}\frac{a}{b}\int_{1}^{R}\frac{d(\frac{b}{s^{\lambda-1}})}{\sqrt{\left(\frac{b}{s^{\lambda-1}}\right)^{2}
-\frac{b^{2}}{R^{2\lambda-2}}}}\right] \\
&=&\left(\frac{1}{R^{\lambda-1}}+\sqrt{\frac{1}{R^{2(\lambda-1)}}-1}\right)^{\frac{1}{1-\lambda}\frac{a}{b}}\\
&=&\left(R^{1-\lambda}+\sqrt{R^{2-2\lambda}-1}\right)^{\frac{1}{1-\lambda}\frac{a}{b}}, \quad \quad(\frac{1}{|w|^{\lambda}}-\text{Nitsche type inequality})
\end{eqnarray*}
\[\eqno{(2.26)}\]
which happens to be (1.16).
Notice that when $\lambda=0$, (2.26) becomes (1.12). For $\alpha\geq-b^{2}$ and $\lambda<1$, we get
\begin{eqnarray*}
&&\sqrt{1+\frac{R^{2(\lambda-1)}\alpha}{b^{2}}}-R^{2(\lambda-1)}\sqrt{1+\frac{\alpha}{b^{2}}} \\
&=&\frac{\left[1-R^{2(\lambda-1)}\right]\left[1+R^{2(\lambda-1)}+\frac{R^{2(\lambda-1)}\alpha}{b^{2}}\right]}{\sqrt{1+\frac{R^{2(\lambda-1)}\alpha}{b^{2}}}
+R^{2(\lambda-1)}\sqrt{1+\frac{\alpha}{b^{2}}}}\\
&>&0.
\end{eqnarray*}
Since $\frac{1}{\lambda-1}<0$, then $\varphi'(\alpha)<0$ (see (2.19)). Thus, $\varphi(\alpha)$ is strictly decreasing. By (2.26) and
\begin{eqnarray*}
\varphi(-b^{2})&=& \left[\frac{R^{\lambda-1}(1+\sqrt{1+\frac{-b^{2}}{b^{2}}})}{1+\sqrt{1+\frac{R^{2\lambda-2}}{b^{2}}(-b^{2})}}
\right]^{\frac{1}{\lambda-1}\frac{a}{b}}\\
&=&\left(R^{1-\lambda}+\sqrt{R^{2-2\lambda}-1}\right)^{\frac{1}{1-\lambda}\frac{a}{b}}.
\end{eqnarray*}
\[\eqno{(2.27)}\]
We obtain there exists a unique $\alpha\geq-b^{2}$, such that $\varphi(\alpha)=r$. Thus we also conclude that (2.23-2.24) from discussion like $\lambda>1$. The inequality (2.16) and (2.26) are unified as (1.16).

  When $\lambda=1$, (2.15) becomes
\[t=\exp\left(\int_{1}^{y}\frac{a}{\sqrt{\alpha+b^{2}}}\frac{1}{s}ds\right)=y^{\frac{a}{\sqrt{\alpha+b^{2}}}}.\eqno{(2.28)}\]
So we get $H(t)=t^{\frac{\sqrt{\alpha+b^{2}}}{a}}$. By $H(r)=R$, we have $\frac{\sqrt{\alpha+b^{2}}}{a}=\frac{\ln R}{\ln r}$. Therefore we obtain the radial mapping
\[h_{1}^{\ast}(te^{i\theta})=t^{\frac{\ln R}{\ln r}}e^{i\theta}.\eqno{(2.29)}\]
And a direct computation yields
\[\mathbb{E}[h_{1}^{\ast}]=2\pi\ln R\left(a^{2}\frac{\ln R}{\ln r}+b^{2}\frac{\ln r}{\ln R}\right).\eqno{(2.30)}\]

\section{Proof of Theorem 1.1}
\subsection{When $\lambda\neq1$}
In view of $J(z,h)=\Im(\overline{h_{_{N}}}h_{T})$ and $|h|_{N}\leq |h_{N}|$,
for $p$, $q$, $A$, $B\in \mathbb{R}$, we obtain the following general inequality
\begin{eqnarray*}
& &\frac{a^{2}|h_{N}|^{2}+b^{2}|h_{T}|^{2}}{|h(z)|^{2\lambda}} \\
&\geq&\frac{(a^{2}-b^{2}p^{2})|h_{N}|^{2}+(b^{2}-a^{2}q^{2})|h_{T}|^{2}+2abpq|h_{N}||h_{T}|}{|h(z)|^{2\lambda}} \\
 &\geq&\frac{2(a^{2}-b^{2}p^{2})At\frac{|h|_{N}}{t}-(a^{2}-b^{2}p^{2})A^{2}+2(b^{2}-a^{2}q^{2})|h|\Im[\frac{h_{T}}{h}]B
 }{|h(z)|^{2\lambda}}\\
 &&\frac{-(b^{2}-a^{2}q^{2})B^{2}+2abpqJ(z,h)  }{|h(z)|^{2\lambda}}.
\end{eqnarray*}
\[\eqno{(3.1)}\]
\textbf{Case 1}: We take $q=\frac{b}{a}$, $p=\frac{\frac{a}{|h|^{\lambda-1}}}{\sqrt{\frac{b^{2}}{|h|^{2(\lambda-1)}}+\alpha}}$,
$A=\frac{\sqrt{\frac{b^{2}}{|h|^{2(\lambda-1)}}+\alpha}}{a t}|h|^{\lambda}$,  and (3.1) becomes
\begin{eqnarray*}
& &\frac{a^{2}|h_{N}|^{2}+b^{2}|h_{T}|^{2}}{|h(z)|^{2\lambda}} \\
    &\geq&\frac{2(a^{2}-b^{2}p^{2})At\frac{|h|_{N}}{t}+2abpqJ(z,h)-(a^{2}-b^{2}p^{2})A^{2}
    }{|h(z)|^{2\lambda}}\\
   &=&\frac{2a\alpha}{\sqrt{\frac{b^{2}}{|h|^{2(\lambda-1)}}+\alpha}}\frac{|h|_{N}}{t}\frac{1}{|h|^{\lambda}}
   +2b^{2}\frac{\frac{a}{|h|^{3\lambda-1}}}{\sqrt{\frac{b^{2}}{|h|^{2(\lambda-1)}}+\alpha}}J(z,h)-\frac{\alpha}{t^{2}}.
\end{eqnarray*}
\[\eqno{(3.2)}\]
Integrating both sides, we obtain
\begin{eqnarray*}
& &\iint_{\mathbb{A}_{1}}\frac{a^{2}|h_{N}|^{2}+b^{2}|h_{T}|^{2}}{|h(z)|^{2\lambda}}dz \\
    &\geq&\iint_{\mathbb{A}_{1}}\frac{2a\alpha}{\sqrt{\frac{b^{2}}{|h|^{2(\lambda-1)}}+\alpha}}\frac{|h|_{N}}{t}\frac{1}{|h|^{\lambda}}dz
   +2b^{2}\iint_{\mathbb{A}_{1}}\frac{\frac{a}{|h|^{3\lambda-1}}J(z,h)}{\sqrt{\frac{b^{2}}{|h|^{2(\lambda-1)}}+\alpha}}dz-\iint_{\mathbb{A}_{1}}\frac{\alpha}{t^{2}}dz\\
    &=&4\pi a\alpha\int_{1}^{r}\frac{\frac{1}{|h|^{\lambda}}}{\sqrt{\frac{b^{2}}{|h|^{2(\lambda-1)}}+\alpha}}\frac{\partial|h|}{\partial t}dt
   +4\pi ab\int_{1}^{R}\frac{\frac{1}{\tau^{3\lambda-2}}}{\sqrt{\frac{1}{\tau^{2(\lambda-1)}}+\frac{\alpha}{b^{2}}}}d\tau-2\pi \alpha\int_{1}^{r}\frac{1}{t}dt\\
   &=&\frac{4\pi a\alpha}{b(1-\lambda)}\ln \frac{\frac{b}{R^{^{\lambda-1}}}+\sqrt{\frac{b^{^{2}}}{R^{2(\lambda-1)}}+\alpha}}{b+\sqrt{b^{2}+\alpha}}+
   \frac{2\pi a\alpha}{b(1-\lambda)}\ln \frac{b+\sqrt{b^{2}+\alpha}}{\frac{b}{R^{\lambda-1}}+\sqrt{\frac{b^{^{2}}}{R^{^{2(\lambda-1)}}}+\alpha}}\\
   & &+ \frac{2\pi a}{(1-\lambda)b}\left[\frac{b}{R^{\lambda-1}}\sqrt{\frac{b^{2}}{R^{2(\lambda-1)}}+\alpha}-b\sqrt{b^{2}+\alpha}\right]-2\pi\alpha\ln r.
\end{eqnarray*}
\[\eqno{(3.3)}\]
It is easy to see from (2.23), (2.25) and (3.3) that
\begin{eqnarray*}
\mathbb{E}_{\lambda}[h]&=&\iint_{\mathbb{A}_{1}}\frac{a^{2}|h_{N}|^{2}+b^{2}|h_{T}|^{2}}{|h(z)|^{2\lambda}}dz \\
    &\geq&\frac{2\pi ab}{\lambda-1}\left(\sqrt{1+\frac{\alpha}{b^{2}}}-\frac{1}{R^{^{\lambda-1}}}\sqrt{\frac{1}{R^{2(\lambda-1)}}+\frac{\alpha}{b^{2}}}\right)\\
   &=&\iint_{\mathbb{A}_{1}}\frac{a^{2}|h^{\ast}_{\lambda N}|^{2}+b^{2}|h^{\ast}_{\lambda T}|^{2}}{|h^{\ast}_{\lambda}(z)|^{2\lambda}}dz\\
   &=&\mathbb{E}_{\lambda}[h^{\ast}_{\lambda}].
\end{eqnarray*}
\[\eqno{(3.4)}\]
$\mathbf{Case}$ $\mathbf{2:}$ We take $p=\frac{a}{b}$, $q=\frac{\sqrt{\frac{b^{2}}{|h|^{2(\lambda-1)}}+\alpha}}{\frac{a}{|h|^{\lambda-1}}}$,
$B=\frac{|h|}{t}$, and (3.1) reduces to
\begin{eqnarray*}
& &\frac{a^{2}|h_{N}|^{2}+b^{2}|h_{T}|^{2}}{|h(z)|^{2\lambda}} \\
    &\geq&\frac{2(b^{2}-a^{2}q^{2})B|h|\Im\frac{h_{T}}{h}+2abpqJ(z,h)-(b^{2}-a^{2}q^{2})B^{2}
    }{|h(z)|^{2\lambda}}\\
    &=&\frac{2\left(b^{2}-a^{2}\frac{\frac{b^{2}}{|h|^{2(\lambda-1)}}+\alpha}{\frac{a^{2}}{|h|^{2(\lambda-1)}}}\right)|h|\frac{|h|}{t}\Im\frac{h_{T}}{h}
    +2a^{2}\frac{\sqrt{\frac{b^{2}}{|h|^{2(\lambda-1)}}+\alpha}}{\frac{a}{|h|^{\lambda-1}}}J(z,h)
    }{|h|^{2\lambda}}\\
    & &\frac{ -\left(b^{2}-a^{2}\frac{\frac{b^{2}}{|h|^{2(\lambda-1)}}+\alpha}{\frac{a^{2}}{|h|^{2(\lambda-1)}}}\right)\frac{|h|^{2}}{t^{2}}}{|h|^{2\lambda}}\\
   &=&-2\alpha\frac{\Im\frac{h_{T}}{h}}{t}+2a\frac{\sqrt{\frac{b^{2}}{|h|^{2(\lambda-1)}}+\alpha}}{|h|^{\lambda+1}}J(z,h)+\frac{\alpha}{t^{2}}.
\end{eqnarray*}
\[\eqno{(3.5)}\]
Integrating this inequality (3.5) over $\mathbb{A}_{1}$, we have
\begin{eqnarray*}
& &\iint_{\mathbb{A}_{1}}\frac{a^{2}|h_{N}|^{2}+b^{2}|h_{T}|^{2}}{|h(z)|^{2\lambda}}dz \\
    &\geq&-2\alpha\iint_{\mathbb{A}_{1}}\frac{\Im\frac{h_{T}}{h}}{t}dz
   +2a\iint_{\mathbb{A}_{1}}\frac{\sqrt{\frac{b^{2}}{|h|^{2(\lambda-1)}}+\alpha}}{|h|^{\lambda+1}}J(z,h)dz+\iint_{\mathbb{A}_{1}}\frac{\alpha}{t^{2}}dz\\
    &=&-4\pi\alpha\int_{1}^{r}\frac{1}{t}dt+\frac{4\pi ab}{1-\lambda}\int_{1}^{R}\sqrt{\frac{1}{\tau^{^{2(\lambda-1)}}}+\frac{\alpha}{b^{2}}}d\left(\frac{1}{\tau^{^{\lambda-1}}}\right)
    +2\pi\alpha\int_{1}^{r}\frac{1}{t}dt\\
    &=&-2\pi\alpha\ln r+\frac{2\pi ab}{1-\lambda}\left(\frac{1}{R^{^{\lambda-1}}}\sqrt{\frac{1}{R^{2(\lambda-1)}}+\frac{\alpha}{b^{2}}}-\sqrt{1+\frac{\alpha}{b^{2}}}\right)\\
   &&+\frac{2\pi a\alpha}{(1-\lambda)b}\ln \frac{\frac{1}{R^{\lambda-1}}+\sqrt{\frac{1}{R^{2(\lambda-1)}}+\frac{\alpha}{b^{2}}}}{1+\sqrt{1+\frac{\alpha}{b^{2}}}}.
\end{eqnarray*}
\[\eqno{(3.6)}\]
By means of (2.23), (2.25) and (3.6), we obtain (3.4) again.

$\mathbf{Uniqueness}$  Let
\[h(te^{i\theta})=\rho(t,\theta)e^{i\varphi(t,\theta)}. \eqno{(3.7)}\]
Firstly, we prove the uniqueness in Case 1. Now suppose the equality in (3.4) holds. Then all the equalities (3.1) and (3.2) must hold, which implies that
\[J(z,h)=\Im (\overline{h_{N}}h_{T})=|h_{N}||h_{T}|, \eqno{(3.8)}\]
\[|h_{N}|=|h|_{N}, \eqno{(3.9)}\]
\[bp|h_{N}|=aq|h_{T}|, \eqno{(3.10)}\]
\[A=|h_{N}|=\frac{\sqrt{\frac{b^{2}}{|h|^{2(\lambda-1)}}+\alpha}}{at}|h|^{\lambda}, \eqno{(3.11)}\]
\[|h_{T}|=|h|\Im[\frac{h_{T}}{h}]. \eqno{(3.12)}\]
By (3.9), we have
\[\rho_{t}=|h|_{N}=|h_{N}|=|h_{t}|=\sqrt{\rho_{t}^{2}+\rho^{2}\varphi_{t}^{2}},\eqno{(3.13)} \]
then $\varphi_{t}=0$. Thus $\varphi(t,\theta)=\varphi(\theta)$. By $q=\frac{b}{a}$, $p=\frac{\frac{a}{|h|^{\lambda-1}}}{\sqrt{\frac{b^{2}}{|h|^{2(\lambda-1)}}+\alpha}}$ and (3.10), we know that $bp|h_{N}|=aq|h_{T}|=b|h_{T}|$, then
\[p|h_{N}|=|h_{T}|.\eqno{(3.14)}\]
By means of (3.11) and (3.14), we get
\[\frac{\frac{a}{|h|^{\lambda-1}}}{\sqrt{\frac{b^{2}}{|h|^{2(\lambda-1)}}+\alpha}}\frac{\sqrt{\frac{b^{2}}{|h|^{2(\lambda-1)}}+\alpha}}{at}|h|^{\lambda}=|h_{T}|\]
which reduces to
\[\frac{|h|}{t}=|h_{T}|.\eqno{(3.15)}\]
Thus
\[\frac{\rho}{t}=|h_{T}|=|\frac{\rho_{\theta}e^{i\varphi}+\rho e^{i\varphi }i\varphi_{\theta}}{t}|=\frac{\sqrt{\rho^{2}_{\theta}+\rho^{2}\varphi^{2}_{\theta}}}{t},\eqno{(3.16)}\]
then
\[\rho=\sqrt{\rho^{2}_{\theta}+\rho^{2}\varphi^{2}_{\theta}}. \eqno{(3.17)}\]
By (3.12) and (3.15-3.17), we obtain
\[\frac{\sqrt{\rho^{2}_{\theta}+\rho^{2}\varphi^{2}_{\theta}}}{t}=\rho\frac{\varphi_{\theta}}{t}
=\sqrt{\rho^{2}_{\theta}+\rho^{2}\varphi^{2}_{\theta}}\frac{\varphi_{\theta}}{t},\]
then
\[\varphi_{\theta}=1.\eqno{(3.18)}\]
According to $\varphi_{t}=0$, we have $\varphi=\theta+\beta$, where $\beta$ is a real constant. By (3.17) and (3.18) we infer that
$\rho=\sqrt{\rho^{2}_{\theta}+\rho^{2}}$, then $\rho_{\theta}=0$. Thus $\rho(t,\theta)=\rho(t)$. By (3.10) and (3.12) we
have
$p|h_{N}|=|h_{T}|=|h|\Im[\frac{h_{T}}{h}],$
and then
$p\rho_{t}=\rho\frac{\varphi_{\theta}}{t}=\frac{\rho}{t},$
that is
$\frac{\frac{a}{|h|^{\lambda-1}}}{\sqrt{\frac{b^{2}}{|h|^{2(\lambda-1)}}+\alpha}}\rho_{t}=\frac{\rho}{t}.$
We have
\[\frac{\frac{a}{(1-\lambda)b}d\left(\frac{b}{\rho^{\lambda-1}}\right)}{\sqrt{\frac{b^{2}}{\rho^{2(\lambda-1)}}+\alpha}}=\frac{dt}{t}.\eqno{(3.19)}\]
After integrating the equality (3.19) and using that $|h|=1$ for $|z|=1$, we arrive at
\[\rho=\frac{2^{\frac{1}{\lambda-1}}(1+\sqrt{1+\frac{\alpha}{b^{2}}})^{\frac{1}{\lambda-1}}t^{\frac{b}{a}}}{\left[(1+\sqrt{1+\frac{\alpha}{b^{2}}})^{2}
-t^{\frac{2(\lambda-1)b}{a}}\frac{\alpha}{b^{2}}\right]^{\frac{1}{\lambda-1}}}. \eqno{(3.20)}\]
Thus
\[h(z)=h^{\ast}_{\lambda}(z)=\frac{2^{\frac{1}{\lambda-1}}(1+\sqrt{1+\frac{\alpha}{b^{2}}})^{\frac{1}{\lambda-1}}|z|^{\frac{b}{a}}}
{\left[(1+\sqrt{1+\frac{\alpha}{b^{2}}})^{2}-|z|^{\frac{2(\lambda-1)b}{a}}\frac{\alpha}{b^{2}}\right]^{\frac{1}{\lambda-1}}}
e^{i\theta+\beta},\eqno{(3.21)}\]
where $\beta$ is a real constant.

The uniqueness in case 2 is obtained in a similar way.

\subsection{When $\lambda=1$}

Obviously, for $p_{\ast}$, $q_{\ast}\in \mathbb{R}$, we have the following inequality
\[
\frac{a^{2}|h_{N}|^{2}+b^{2}|h_{T}|^{2}}{|h(z)|^{2}}\geq\frac{(a^{2}-b^{2}p_{\ast}^{2})|h_{N}|^{2}+(b^{2}-a^{2}q_{\ast}^{2})|h_{T}|^{2}+2p_{\ast}q_{\ast}ab|h_{N}|h_{T}|}
{|h(z)|^{2}}.\eqno{(3.22)}\]
We also consider two cases.

\textbf{Case 1:} Let $q_{\ast}=\frac{b}{a}$, and (3.22) reduces to
\[
\frac{a^{2}|h_{N}|^{2}+b^{2}|h_{T}|^{2}}{|h(z)|^{2}}\geq\frac{(a^{2}-b^{2}p_{\ast}^{2})|h_{N}|^{2}+2p_{\ast}b^{2}|h_{N}||h_{T}|}
{|h(z)|^{2}}.\eqno{(3.23)}
\]
Since $|h_{N}||h_{T}|\geq\Im(\overline{h_{N}}h_{T})=J(z,h)$
and
$|h_{N}|\geq|h|_{N},$
 by (3.23) we have
\[
\frac{a^{2}|h_{N}|^{2}+b^{2}|h_{T}|^{2}}{|h(z)|^{2}}\geq\frac{(a^{2}-b^{2}p_{\ast}^{2})|h|_{N}^{2}
+2p_{\ast}b^{2}J(z,h)}{|h(z)|^{2}}.\eqno{(3.24)}
\]
By integrating over $\mathbb{A}_{1}$, we can get
\[
\iint_{\mathbb{A}_{1}}\frac{a^{2}|h_{N}|^{2}+b^{2}|h_{T}|^{2}}{|h(z)|^{2}}dz\geq(a^{2}-b^{2}p_{\ast}^{2})
\iint_{\mathbb{A}_{1}}\frac{|h|_{N}^{2}}{|h(z)|^{2}}dz+2p_{\ast}b^{2}\iint_{\mathbb{A}_{1}}\frac{J(z,h)}{|h(z)|^{2}}dz.
\eqno{(3.25)}\]
By H$\ddot{{\rm o}}$lder inequality, we have
\begin{eqnarray*}
\Big(\iint_{\mathbb{A}_{1}}\frac{|h|_{N}}{|z||h(z)|}dz\Big)^{2}&\leq&\iint_{\mathbb{A}_{1}}\frac{|h|_{N}^{2}}{|h(z)|^{2}}dz
\iint_{\mathbb{A}_{1}}\frac{1}{|z|^{2}}dz\\
   &=&2\pi \ln r \iint_{\mathbb{A}_{1}}\frac{|h|_{N}^{2}}{|h(z)|^{2}}dz.
\end{eqnarray*}
\[\eqno{(3.26)}\]
By Fubini's theorem
\begin{eqnarray*}
\iint_{\mathbb{A}_{1}}\frac{|h|_{N}}{|z||h(z)|}dz&=&\int^{2\pi}_{0}\left(\int_{1}^{r}\frac{1}{|h(z)|}\frac{\partial|h(z)|}{\partial t}dt\right)d\theta\\
&=&\int^{2\pi}_{0}\left(\int_{1}^{R}\frac{d\varrho}{\varrho}\right)d\theta=2\pi \ln R.
\end{eqnarray*}
\[\eqno{(3.27)}\]
According to (3.26) and (3.27), we know
\[
\iint_{\mathbb{A}_{1}}\frac{|h|_{N}^{2}}{|h(z)|^{2}}dz\geq 2\pi\frac{\ln^{2}R}{\ln r}.\eqno{(3.28)}
\]
Moreover,
\[
\iint_{\mathbb{A}_{1}}\frac{J(z,h)}{|h(z)|^{2}}dz=\iint_{\mathbb{A}_{2}}\frac{1}{|w|^{2}}dw=2\pi\ln R.\eqno{(3.29)}
\]
According to (3.25), (3.28) and (3.29), we obtain
\[
\iint_{\mathbb{A}_{1}}\frac{a^{2}|h_{N}|^{2}+b^{2}|h_{T}|^{2}}{|h(z)|^{2}}dz\geq2\pi(a^{2}-b^{2}p_{\ast}^{2})\frac{\ln^{2} R}{\ln r}+ 2p_{\ast}b^{2}2\pi\ln R.
\eqno{(3.30)}\]
Let $p_{\ast}=\frac{\ln r}{\ln R}$, and by (2.30), then (3.30) reduces to
\begin{eqnarray*}
\mathbb{E}_{1}[h]&=&
\iint_{\mathbb{A}_{1}}\frac{a^{2}|h_{N}|^{2}+b^{2}|h_{T}|^{2}}{|h(z)|^{2}}dz\\
&\geq&2\pi\ln R\Big(a^{2}\frac{\ln R}{\ln r}+b^{2}\frac{\ln r}{\ln R}\Big)\\
&=&\iint_{\mathbb{A}_{1}}\frac{a^{2}|h^{*}_{N}|^{2}+b^{2}|h^{*}_{T}|^{2}}{|h_{1}^{*}(z)|^{2}}dz\\
&=&\mathbb{E}_{1}[h^{\ast}_{1}].
\end{eqnarray*}
\textbf{Case 2:} We take $p_{\ast}=\frac{a}{b}$ in (3.22) and obtain
\begin{eqnarray*}
\frac{a^{2}|h_{N}|^{2}+b^{2}|h_{T}|^{2}}{|h(z)|^{2}}
&\geq&\frac{(b^{2}-a^{2}q_{\ast}^{2})|h_{T}|^{2}+2a^{2}q_{\ast}|h_{N}||h_{T}|}{|h(z)|^{2}}\\
&\geq&\frac{(b^{2}-a^{2}q_{\ast}^{2})|h_{T}|^{2}+2a^{2}q_{\ast}J(z,h)}{|h(z)|^{2}}.
\end{eqnarray*}
\[\eqno{(3.31)}\]
Integrating over $\mathbb{A}_{1}$ in (3.31), we get
\[
\iint_{\mathbb{A}_{1}}\frac{a^{2}|h_{N}|^{2}+b^{2}|h_{T}|^{2}}{|h(z)|^{2}}dz
\geq(b^{2}-a^{2}q_{\ast}^{2})\iint_{\mathbb{A}_{1}}\frac{|h_{T}|^{2}}{|h(z)|^{2}}dz+2a^{2}q_{\ast}\iint_{\mathbb{A}_{1}}\frac{J(z,h)}{|h(z)|^{2}}dz.\eqno{(3.32)}
\]
By means of H$\ddot{{\rm o}}$lder inequality, we have
\begin{eqnarray*}
\Big(\iint_{\mathbb{A}_{1}}\frac{|h_{\theta}|}{|z|^{2}|h(z)|}dz\Big)^{2}
&\leq&\iint_{\mathbb{A}_{1}}\frac{|h_{\theta}|^{2}}{|z|^{2}|h(z)|^{2}}dz\iint_{\mathbb{A}_{1}}\frac{1}{|z|^{2}}dz\\
&=&2\pi \ln r\iint_{\mathbb{A}_{1}}\frac{|h_{\theta}|^{2}}{|z|^{2}|h(z)|^{2}}dz.
\end{eqnarray*}
\[\eqno{(3.33)}\]
Since
\[
2\pi\ln r=2\pi \int_{1}^{r}\frac{dt}{t}=\iint_{\mathbb{A}_{1}}\frac{1}{|z|}\Im\Big(\frac{h_{T}}{h}\Big)dz
\leq\iint_{\mathbb{A}_{1}}\frac{|h_{\theta}|}{|z|^{2}|h(z)|}dz.\eqno{(3.34)}
\]
According to (3.33) and (3.34), we get
\[
\iint_{\mathbb{A}_{1}}\frac{|h_{T}|^{2}}{|h(z)|^{2}}dz=\iint_{\mathbb{A}_{1}}\frac{|h_{\theta}|^{2}}{|z|^{2}|h(z)|^{2}}dz
\geq2\pi\ln r.\eqno{(3.35)}
\]
Therefore, (3.32) reduces to
\[
\iint_{\mathbb{A}_{1}}\frac{a^{2}|h_{N}|^{2}+b^{2}|h_{T}|^{2}}{|h(z)|^{2}}dz
\geq(b^{2}-a^{2}q_{\ast}^{2})2\pi \ln r+2a^{2}q_{\ast}2\pi\ln R.\eqno{(3.36)}
\]
Let $q_{\ast}=\frac{\ln R}{\ln r}$, by (3.36) we get
\begin{eqnarray*}
\mathbb{E}_{1}[h]&=&\iint_{\mathbb{A}_{1}}\frac{a^{2}|h_{N}|^{2}+b^{2}|h_{T}|^{2}}{|h(z)|^{2}}dz\\
&\geq&2\pi\ln R\Big(a^{2}\frac{\ln R}{\ln r}+b^{2}\frac{\ln r}{\ln R}\Big)\\
&\geq&\iint_{\mathbb{A}_{1}}\frac{a^{2}|h^{*}_{N}|^{2}+b^{2}|h^{*}_{T}|^{2}}{|h_{1}^{*}(z)|^{2}}dz\\
&=&\mathbb{E}_{1}[h^{\ast}_{1}].
\end{eqnarray*}
As in the proof of $\lambda\neq1$, the uniqueness is obtained in a similar way.
The proof of Theorem 1.1 is completed.

\section{The extremal problem for weighted combined distortion}

In this section, we will consider the extremal problem for weighted combined distortion. Certainly, the extremal problem for integrable distortion  is connected with the extremal problem for harmonic energy.  For the related approaches to the extremal problem for integrable distortion we refer to the papers (\cite{AIMO},\cite{AIM}, \cite{FTWS}, \cite{IMO}, \cite{MM} etc.).

At first, we give the definition of the weighted combined distortion.
For $a,b>0$ and $h\in\mathfrak{F}(\mathbb{A}_{1},\mathbb{A}_{2})$, we define the weighted combined distortion by
 \[\mathbb{K}_{\lambda}[h]=\iint_{\mathbb{A}_{1}}\frac{a^{2}\rho^{2}|\nabla\Theta|^{2}
 +b^{2}|\nabla\rho|^{2}}{J(z,h)|z|^{2\lambda}}dz.\eqno{(4.1)}\]
where $h(z)=\rho e^{i\Theta}$.

Next, by similar reasoning of in proof Theorem 4.1 in \cite{YTF}, we obtain the relationship between weighted combined energy and weighted combined distortion.
\begin{lemma}
Assume that $h$ is a map of $\mathbb{A}_{1}$ to $\mathbb{A}_{2}$, and let $f=h^{-1}$. If
\[E_{\lambda}[f]=\iint_{\mathbb{A}_{2}}(a^{2}|f_{N}|^{2}+b^{2}|f_{T}|^{2})\frac{1}{|f(w)|^{2\lambda}}dw,\eqno{(4.2)}\]
then
\[E_{\lambda}[f]=\mathbb{K}_{\lambda}[h].\eqno{(4.3)}\]
where $\mathbb{K}_{\lambda}[h]$ is defined in (4.1).
\end{lemma}

Finally, we can obtain the following result from Theorem 1.1 and Lemma 4.1.

\begin{theorem}
(I) When $\lambda\neq1$, under the Nitsche type condition (1.16), for $f\in\mathfrak{H}(\mathbb{A}_{2},\mathbb{A}_{1})$,
  the weighted combined distortion $\mathbb{K}_{\lambda}[f]$ attains its minimum for a radial mapping $f^{\ast}_{\lambda}=h_{\lambda}^{\ast-1}$ (as defined (1.17)). That is to say,
\[\mathbb{K}_{\lambda}[f]\geq\mathbb{K}[f^{\ast}_{\lambda}]. \]
The minimizer is unique up to a rotation of the annuli.

(II) When $\lambda=1$, for $f\in\mathfrak{H}(\mathbb{A}_{2},\mathbb{A}_{1})$,
  the weighted combined distortion $\mathbb{K}_{\lambda}[f]$ attains its minimum for a radial mapping
  $f^{\ast}_{1}=h_{1}^{\ast-1}$ (as defined (1.20)).
The minimizer is unique up to a rotation of the annuli.
\end{theorem}

 \bigskip

{\bf Acknowledgements} \quad  The authors would like to thank to the referee for a very careful reading of manuscript.

\bibliographystyle{amsplain}

\end{document}